\documentclass[12pt]{article}
\usepackage{mathtools,amsfonts,amssymb,amsthm,amscd}
\usepackage{microtype}
\usepackage{color}
\usepackage[height=9.3in,width=6.5in]{geometry}

\usepackage{tikz}
\usepackage{latexsym}
\usepackage[pdftitle={Sandpiles and unicycles on random planar maps},
            pdfauthor={Xin Sun and David B. Wilson},
            bookmarks=true,bookmarksopen=true,bookmarksopenlevel=3]{hyperref}

\newtheorem{theorem}{Theorem}[section]

\newtheorem{lemma}[theorem]{Lemma}
\newtheorem{proposition}[theorem]{Proposition}

\newtheorem{remark}[theorem]{Remark}

\newcommand{\tquad}{t_{\textup{\textrm{quad}}}}
\newcommand{\squad}{s_{\textup{\textrm{quad}}}}

\newcommand*\bcircled[1]{\tikz[baseline=(char.base)]{
            \node[shape=circle,draw,inner sep=0.5pt] (char) {#1};}}
\newcommand*\squared[1]{\tikz[baseline=(char.base)]{
            \node[shape=rectangle,draw,inner sep=2pt] (char) {#1};}}
\def\frshorder{\squared{\textup{\textsf{\hspace{-.1ex}F\hspace{.2ex}}}}}
\def\hambg{\bcircled{\textup{\textsf{H}}}}
\def\hamdr{\squared{\textup{\textsf{H}}}}
\def\chsdr{\squared{\textup{\textsf{\hspace{-.3ex}C\hspace{.3ex}}}}}
\def\chsbg{\bcircled{\textup{\textsf{\hspace{-.3ex}C\hspace{.3ex}}}}}
\newcommand{\be}{\begin{equation}}
\newcommand{\ee}{\end{equation}}
\def\ba{\begin{align}}
\def\ea{\end{align}}
\newcommand{\G}{\mathcal{G}}
\newcommand{\Rec}{\mathcal{R}}
\newcommand{\T}{\mathcal{T}}
\newcommand{\U}{\mathcal{U}}

\newcommand{\M}{\mathcal{M}}

\renewcommand{\L}{\mathcal{L}}
\newcommand{\1}{\mathbf{1}}

\newcommand{\E}{\mathbb{E}}

\newcommand{\Z}{\mathbb{Z}}

\renewcommand{\P}{\mathbb{P}}

\newcommand{\eps}{\varepsilon}
\def\P{\mathbb{P}}
\def\E{\mathbb{E}}

\DeclareMathOperator{\Var}{Var}
\DeclareMathOperator{\level}{level}

\begin {document}
\title{Sandpiles and unicycles on random planar maps}
\author{\begin{tabular}{c}\href{http://math.mit.edu/~xinsun89}{Xin Sun}\\[-4pt]\small Massachusetts Institute of Technology\end{tabular} \and \begin{tabular}{c}\href{http://dbwilson.com}{David B. Wilson}\\[-4pt]\small Microsoft Research\end{tabular}}
\date{}
\maketitle

\begin{abstract}
We consider the abelian sandpile model and the uniform spanning unicycle on random planar maps.  We show that the sandpile density converges to 5/2 as the maps get large.  For the spanning unicycle, we show that the length and area of the cycle converges to the hitting time and location of a simple random walk in the first quadrant.  The calculations use the ``hamburger-cheeseburger'' construction of Fortuin--Kasteleyn random cluster configurations on random planar maps.
\begin{flushleft}
\end{flushleft}
\end{abstract}

\section{Introduction} \label{sec::The model}
Random planar maps together with discrete statistical mechanics models (e.g., spanning tree, Ising model, FK model) on them is an active research area
 (see e.g., \cite{A03,LM:12,DS:11,Ben13,sheffield2011quantum,GJSZ:12,curien-kortchemski,Wedge,burger-cone,burger-exponent,chen-fk,burger-cone2,burger-cone3}).
In the annealed distribution of a discrete model on a random planar map of a given class,
the joint distribution of the pair $(M,\Sigma)$, where $M$ is the planar map and $\Sigma$ is a configuration of the discrete model on $M$,
is just proportional to the weight of $\Sigma$.
Equivalently, the random map in the class is sampled according to the \textit{partition function\/} of the discrete model (i.e., the weighted sum of all the configurations on the given map), and then a configuration is sampled according to the weighting rule of the discrete model on the map.

We consider two discrete models on random planar maps: the uniform recurrent sandpile and the uniform spanning unicycle (also known as a cycle-rooted spanning tree).  For the sandpile model we show that the sandpile density converges to $5/2$ and concentrates around this value.  For the unicycle model, we compute the weak limit of the joint distribution of the length of the cycle and the area inside the cycle.  The sandpile model calculations depend on the results for the spanning unicycle.

\subsection{Planar maps, spanning trees, and unicycles}\label{sec::intro of maps}
Planar graphs are graphs that can be embedded into the sphere.  We allow the existence of self loops and multiple edges.  A \textit{planar map\/} is a connected planar graph embedded on the sphere considered up to isotopic deformation of the edges, i.e., a planar map contains only information about the combinatorial structure of the embedding.  A \textit{rooted planar map\/} $(M,e)$ is a planar map $M$ with a distinguished directed edge $e$.  We let $\M_n$ denote the set of rooted planar maps with $n$ edges.

Given a finite connected graph $\G=(V,E)$, a \textit{spanning forest\/} of $\G$ is  a subgraph whose vertex set is $V$ and which contains no cycles.  A spanning forest that contains $k$ connected components is called a \textit{$k$-component spanning forest}. A 1-component spanning forest is called a \textit{spanning tree}.  A \textit{$k$-excess subgraph\/} of $\G$ is the union of a spanning tree and $k$ extra edges in $E$. A $1$-excess subgraph is called a \textit{spanning unicycle}.  The planar dual of a $k$-excess subgraph is a $(k+1)$-component spanning forest.
Let $\T(\G)$ denote
\[
\T(\G) := \big\{\text{spanning trees of $\G$}\big\}\,.
\]
and $\U_k(\G)$ denote
\[
\U_k(\G) := \big\{\text{$k$-excess subgraphs of $\G$}\big\}\,.
\]

Given a rooted planar map $(M,e)$, since $e$ is a directed edge, there is a unique face $f$ to the right of $e$, which we call the \textit{outer face}.
The outer face allows us to distinguish between the two sides of a cycle: the \textit{outside\/} is the side containing the outer face, and the \textit{inside\/} is the side which does not.
The \textit{length\/} of a cycle is the number of edges on the cycle.  We define the \textit{area\/} of a cycle to be twice the number of edges inside the cycle plus its length.
Planar maps are in natural bijective correspondence with quadrangulations, and
this combinatorial definition of area corresponds to assigning each quadrangle area~$2$.

\subsection{The abelian sandpile model}
The abelian sandpile model is a model for self-organized criticality \cite{BTK} which is defined as follows.  (See also \cite{HLMPPW} for further background.)
Suppose $\G=(V,E)$ is a finite connected undirected graph with loops and multiple edges allowed.  Let $c(v,w)$ be the number of edges between vertices $v$ and $w$, where self-loops count twice. For $v\in V$, the degree of $v$ is denoted by $\deg(v)=\sum_{w\in V}c(v,w)$. A sandpile configuration on the graph $\G$, with respect to a distinguished vertex $s$ called the sink, assigns a non-negative integer number of grains of sand to each vertex other than the sink $s$.  If a vertex $v\neq s$ has more sand than its degree, then $v$ is \textit{unstable}, and may \textit{topple}, sending one grain of sand to each neighbor.  The sink $s$ never topples.  Since every vertex is connected to the sink, we may repeatedly topple unstable vertices until every vertex is stable. The resulting sandpile is called the
stabilization of the original sandpile, and is independent of the order in which vertices are toppled (which is the abelian property).

Some sandpile configurations are \textit{recurrent}, meaning that from any sandpile configuration, it is possible to add some
amount of sand to the vertices and stabilize to obtain the given configuration. These sandpile configurations are the recurrent states of the Markov chain which at each step adds a grain of sand to a random vertex and then stabilizes the configuration.
The stationary distribution of this Markov chain is the uniform distribution on recurrent sandpile configurations.

We let $\Rec(\G,s)$ denote the set of the recurrent sandpile configurations on a graph $\G$ with sink~$s$.  Majumdar and Dhar gave a bijection between $\Rec(\G,s)$ and $\T(\G)$ \cite{majumdar-dhar:height}. In particular, $|\Rec(\G,s)|=|\T(\G)|$.

Given a recurrent sandpile configuration $\sigma \in \Rec(\G,s)$, where $\sigma(v)$ for $v\neq s$ is the number of grains at $v$,
it is convenient to define $\sigma(s)=\deg(s)$.
Then the total amount of sand (i.e., including the sand at the sink) is
\[
|\sigma| = \sum_{v\in\G} \sigma(v)\,.
\]
With this convention, the distribution of $|\sigma|$ for a random recurrent sandpile $\sigma\in\Rec(\G,s)$
does not depend on the choice of the sink~$s$ \cite{merino-lopez}.
This distribution can be understood in terms of the $k$-excess subgraphs of $\G$ \cite{merino-lopez},
which we will explain in Section~\ref{sec:Leveland Tutte polynomial}.

We define the \textit{sandpile edge density\/} for~$\G$ to be
\begin{equation}\label{eq: edge density}
\rho_e(\G)=\frac{1}{|E|\cdot|R(\G,s)|} \sum\limits_{\sigma\in \Rec(\G,s) } |\sigma|
\end{equation}
and the \textit{sandpile vertex density} for $\G$ to be
\begin{equation}\label{eq: vertex density}
\rho_v(\G)=\frac{1}{|V|\cdot|R(\G,s)|} \sum\limits_{\sigma\in \Rec(\G,s) } |\sigma|\,.
\end{equation}
These sandpile densities $\rho_e(\G)$ and $\rho_v(\G)$ are independent of the choice of the sink.

\subsection{Sandpiles and $k$-excess subgraphs on random planar maps}\label{sec:models and results}
We let $\P_n^k$ denote the annealed distribution on $k$-excess subgraphs of a rooted planar map containing $n$ edges,
i.e., a random planar map $(M,e)$ with a $k$-excess subgraph $U_k$ on it are chosen with probability
\begin{align}\label{eq:def of measure}
\P_n^k[M,e,U_k] &=\frac{1}{\sum\limits_{(M',e')\in \M_n} |\U_k(M')|}\,.
\end{align}
We let $\E^k_n$ and $\Var^k_n$ denote the expectation and variance with respect to $\P_n^k$.

For the directed edge $e$, we let $\underline e$ denote its source vertex, and $\overline e$ denote its destination vertex.
Because of the bijection between recurrent sandpile configurations $\Rec(M,\underline e)$ and spanning trees $\T(M)=\U_0(M)$,
we can interpret $P_n^0$ as being the annealed distribution of uniform recurrent sandpiles on $\M_n$ with sink at the source of the root edge.

\begin{theorem}\label{thm::sandpile-density}
For the uniform recurrent sandpile on a random planar map with $n$ edges, the sandpile density satisfies
\begin{equation}\label{eq:denisty result}
\lim\limits_{n\to\infty }\E_n^0[\rho_v]=2\lim\limits_{n\to\infty }\E_n^0[\rho_e]=\frac{5}{2},\quad\quad\quad\quad \lim\limits_{n\to\infty} \Var^0_n[\rho_v]=\lim\limits_{n\to\infty} \Var^0_n[\rho_e]=0.
\end{equation}
\end{theorem}	
This sandpile density computation can be compared to the uniform recurrent sandpile density on $\Z^2$.
	In 1994 Grassberger conjectured that the (per vertex) sandpile density on $\Z^2$ is $17/8$, based on the numerical integration of singular 4-dimensional integral expressions given by Priezzhev \cite{priezzhev:heights} for the sandpile height distribution at a vertex.  These integral expressions were greatly simplified by Jeng, Piroux, and Ruelle \cite{JPR}, who verified by numerical integration that the sandpile density for $\Z^2$ is $17/8 \pm 10^{-12}$.  An alternative formulation of the sandpile density, relating it to spanning unicyclic graphs and loop-erased random walk, was given by Poghosyan and Priezzhev \cite{PP:54} and Levine and Peres \cite{LP:54}, which enabled its rigorous exact evaluation in \cite{PPR:54,KW4}.  The density of $17/8$ corresponds to an edge density of $1+\frac{1}{16}$, versus $1+\frac{1}{4}$ for the sandpile on a random planar map.  More recently, a simpler proof of the sandpile density on $\Z^2$ was given by Kassel and Wilson \cite{kassel2014looping}, who then computed the sandpile density for numerous other lattices as well.
	
\bigskip
Our results for unicycles, and more generally $k$-excess subgraphs, can be expressed
in terms of simple random walk in the first quadrant of $\Z^2$.
Suppose $(X_t,Y_t)$ is a simple random walk on $\Z^2$ started from $(1,1)$. Denote
\begin{align}\label{eq:cone-quantity}
\tquad=\inf\{t: X_tY_t=0 \},\quad\quad\quad\quad \squad=X_{\tquad}+Y_{\tquad}\,.
\end{align}
In words, $\tquad$ is the exit time from the first quadrant,
and $\squad$ is the distance from the origin when the walk exits the quadrant.
Random walk in the first quadrant has been studied quite extensively (see e.g., \cite{QuadrantBook,QuadrantUnified, RWCone,MartinBoundary});
for an introductory account see \cite[Section 8.1.3]{lawler-limic}.
Here we record that 
\begin{equation}\label{eq:asymptotic}
\P[\tquad>j]\sim \frac{4}{\pi j}   \quad\quad\mathrm{and}\quad\quad
\P[\squad>\ell]\sim \frac{4}{\pi \ell^2}
\end{equation}
asymptotically as $j\to\infty$ and $\ell\to\infty$ (see \cite[Sec.~3]{FR} and \cite[Example 3, Sec.~1.3]{RWCone}).

\begin{theorem}\label{thm::loop length result}\label{thm:k-loops}
Consider the uniform $k$-excess subgraph of a random planar map under the distribution of $\P_n^k$.
As $n\to\infty$ with $k$ fixed, with probability $1-o_n(1)$ the $k$-excess subgraph has $k$ unnested loops.
Let $L_1,A_1,\dots,L_k,A_k$ denote the length and area of the $k$ loops.
Then the joint distribution of $\{(L_i,A_i)\}_{1\le i\le k}$ converges weakly to $k$ i.i.d.\ samples of the random variables
$(\squad,\tquad)$ described above.
\end{theorem}

The area distribution can be compared to the area of the cycle of a uniform spanning unicycle in $\Z^2$.
The moments of the area were computed for $n\times n$ boxes in $\Z^2$ \cite{KKW}, and these moments suggest that $\P[A>j]\asymp 1/j$.

\bigskip
To prove these theorems, we start in Section~\ref{sec:Leveland Tutte polynomial} by explaining how the $k^{\text{th}}$ moments of the amount of sand in a uniform recurrent sandpile are related to $k$-excess subgraphs.
Then in Section~\ref{sec:burger bijection} we review the ``hamburger-cheeseburger'' bijection,
which constructs rooted planar maps with $n$ edges together with a Fortuin--Kasteleyn configuration.
(These FK configurations are more general than $k$-excess subgraphs.)
This bijection is due to Mullin \cite{Mullin} in the special case where the FK model is a spanning tree,
and to Bernardi and Sheffield \cite{Bernardi,sheffield2011quantum} in general.
We use the formulation in \cite{sheffield2011quantum} since it is more convenient for our purposes.
In Section~\ref{sec:edge-density}, we do some asymptotic analysis of the hamburger-cheeseburger bijection to prove Theorem~\ref{thm::loop length result}.
We conclude the proof of Theorem~\ref{thm::sandpile-density} in Section~\ref{sec:small terms}.

\bigskip
\noindent
\textbf{Acknowledgments:}
We thank Adrien Kassel for helpful comments. This work was begun while the first author was an intern in the theory group at Microsoft Research Redmond,
and was completed at the Newton Institute.

\section{Sandpiles, spanning trees, and the Tutte polynomial}
\label{sec:Leveland Tutte polynomial}
One approach to computing the sandpile density is via the Tutte polynomial.
The Tutte polynomial of an undirected graph $\G=(V,E)$ is a polynomial in two variables defined by
\begin{equation}\label{eq:Tutte polynomial}
T_\G(x,y)=\sum\limits_{E'\subset E} (x-1)^{\kappa(E')-\kappa(E)}\,(y-1)^{\kappa(E')+|E'|-|V|}
\end{equation}
where $\kappa(E')$ is the number of connected components of the spanning subgraph of $\G$
with edge set $E'$.

For the abelian sandpile model on a finite connected graph $\G=(V,E)$ with sink $s$, Biggs defined the \textit{level\/} of a sandpile configuration to be
\[\level(\sigma)=|\sigma|-|E|\]
and showed that $0 \le \level(\sigma) \le |E|-|V|+1$ and that these bounds are tight. Consequently, the sandpile edge density satisfies \[1\le \rho_e\le 2\,.\]
Biggs conjectured and Merino proved \cite{merino-lopez} that the generating function of recurrent sandpiles by level is
the Tutte polynomial evaluated at $(1,y)$:
\begin{equation}\label{eq:level generating function}
\sum\limits_{\sigma\in \Rec(\G,s)} y^{\level(\sigma)} =T_\G(1,y)\,.
\end{equation}
Notice that the generating function is independent of the choice of the sink~$s$.

Since we are assuming that $\G$ is connected,
from~\eqref{eq:level generating function} and \eqref{eq:Tutte polynomial} we see
\begin{equation}\label{eq:level-sum}
  \sum_{\text{sandpiles $\sigma$ of $\G$}} \binom{\level(\sigma)}{\ell} = |\U_\ell(\G)|\,,
\end{equation}
so for a random recurrent sandpile~$\sigma\in\Rec(\G,s)$,
\begin{equation}\label{eq:level-average}
  \E^\G\left[\binom{\level(\sigma)}{\ell}\right] = \frac{|\U_\ell(\G)|}{|\T(\G)|}\,.
\end{equation}
When the graph $\G$ is itself random, in particular a random planar map on $n$ edges, we have
\begin{equation}\label{eq:ell-cumulant}
  \E_n^0\left[\binom{\level(\sigma)}{\ell}\right] = \E^0_n \left[\frac{|\U_\ell(\M_n)|}{|\T(\M_n)|}\right]\,.
\end{equation}
This equation relates the binomial moments of the sandpile level to quantities that can be evaluated by the hamburger-cheeseburger bijection which we will describe in Section~\ref{sec:burger bijection}.

For any random variable $Z$, the moments $\E[Z^\ell]$ of $Z$ can be expressed as linear combinations of the binomial moments $\E[\binom{Z}{\ell}]$,
and the cumulants $\E[(Z-\E[Z])^\ell]$ can be expressed as linear combinations of products of the form $\E[Z^m]\E[Z]^{\ell-m}$.
The variance in particular is $\Var[Z] = 2\,\E[\binom{Z}{2}] + \E[\binom{Z}{1}] - \E[\binom{Z}{1}]^2$.
For a random recurrent sandpile $\sigma$ on a random map,
\begin{align}
\E_n^0[\level(\sigma(M_n))]&=\E_n^0\left[\frac{|\U_1(M_n)|}{|\T(M_n)|}\right] \label{eq:E[sand]}\\
\Var_n^0[\level(\sigma(M_n))]&=2\,\E_n^0\left[\frac{|\U_2(M_n)|}{|\T(M_n)|}\right] + \E_n^0\left[\frac{|\U_1(M_n)|}{|\T(M_n)|}\right]
-\left(\E_n^0\left[\frac{|\U_1(M_n)|}{|\T(M_n)|}\right]\right)^2\,. \label{eq:Var[sand]}
\end{align}
Next we evaluate the terms on the right-hand side using the hamburger-cheesburger bijection.

\section{The hamburger-cheeseburger bijection}\label{sec:burger bijection}

Sheffield \cite[Section 4.1]{sheffield2011quantum} constructed a
bijection called the \textit{hamburger-cheeseburger bijection\/}
between ``perfect words'' over a five-letter alphabet
$\{\,\hambg,\hamdr,\chsbg,\chsdr,\frshorder\,\}$ and FK configurations
on rooted planar maps.
The letters in this alphabet can be interpreted as events at a burger restaurant.
$\hambg$ and $\chsbg$ indicate that a new hamburger or cheeseburger is produced.
New burgers are placed on top of a burger stack, which is initially empty.
$\hamdr$ and $\chsdr$ indicate that a customer has ordered a hamburger or cheeseburger,
in which case the topmost burger of the appropriate type is removed from the stack and
given to the customer.  $\frshorder$ indicates a fresh order, where the customer orders
whichever burger is on top of the stack, regardless of type.
A \textit{perfect word\/} is a sequence of these letters for which every burger order
is fulfilled by a burger already on the stack, and for which the burger stack ends empty.
A perfect word of order~$n$ is one which contains $n$~burgers and $n$~orders, i.e., has
length~$2n$.

The hamburger-cheeseburger bijection maps a perfect word of order $n$ to a rooted planar map
$(M_n,e)$ with $n$ edges together with a subset $E'$ of the edges in the map.
(The edge subset $E'$ is the FK configuration on the map.)
In addition to Sheffield's description of the bijection \cite[Section 4.1]{sheffield2011quantum},
a nice exposition is given by Chen \cite{chen-fk}.
We note here a few basic properties of the bijection:
\begin{enumerate}
\item The number of edges in $E'$ is the number of $\hamdr$'s plus the number of $\frshorder$'s matching $\chsbg$'s.
\item The number of connected components in the subgraph spanned by $E'$ is 1 plus the number of $\frshorder$'s matching $\hambg$'s.
\item Let $(E')^*$ denote the dual edges of $E\setminus E'$ on the dual map $M^*_n$ of $M_n$.
The number of connected components in the dual subgraph spanned by $(E')^*$ is 1 plus the number of $\frshorder$'s matching $\chsbg$'s.
\item
Suppose that a $\frshorder$ matches a $\chsbg$ in a perfect word,
and that in between there are $\ell$ $\hamdr$'s which are fulfilled by \hambg's before the fresh $\chsbg$, and $2m$ other letters.  Then the $\frshorder$
corresponds the ``last edge'' of a loop of $E'$ which has length $\ell+1$,
the portion of the map $M_n$ inside the loop has area
$(1+\ell+2m)$, and the portion of $M_n$ and $E'$ inside of the loop are determined by the subword between the $\frshorder$ and its matching $\chsbg$.
\end{enumerate}


\begin{proposition}\label{prop:bijection}
Let
\begin{equation}
\Theta_n^k = \big\{\text{perfect words of order $n$ with exactly $k$ \frshorder's, which are all fulfilled by \chsbg's}\big\}\,.
\end{equation}
  Under the hamburger-cheeseburger bijection, elements in $\Theta^k_n$
  correspond to triples $(M_n,e, U_k)$, where $(M_n,e)\in \M_n$ and
  $U_k$ is a $k$-excess connected subgraph of $M_n$.
  Furthermore,
\begin{equation} \label{eq:Theta}
\E_n^0\left[\frac{|\U_k(M_n)|}{|\T(M_n)|}\right] = \frac{|\Theta^k_n|}{|\Theta^0_n|}
\end{equation}
\end{proposition}
\begin{proof}
Suppose $E'$ has 1 component and its dual $(E')^*$ has $k+1$ components.  Then $(E')^*$ contains no cycles, so it is a $(k+1)$-component spanning forest, and consequently $E'$ is a $k$-excess subgraph of $M_n$.
%

Recall that under the distribution $\P_n^0$, each rooted map $(M_n,e)$ occurs with probability $|\T(M_n)|/\sum_{(M'_n,e')\in\M_n} |\T(M'_n)| = |\T(M_n)|/|\Theta^0_n|$.  In the expectation in \eqref{eq:Theta}, the $|\T(M_n)|$ terms cancel, giving the right-hand side of \eqref{eq:Theta}.
\end{proof}

\newcommand{\old}[1]{}
\old{
\subsection{Associated quadrangulation and triangulations}\label{sec:quad}

Let $V = V(M)$ be the set of vertices of a planar map $M$ and $F = F(M)$ the set of faces. Let $Q = Q(M)$ be the associated \textit{quadrangulation}, which is the map whose vertex set is $V \cup F$, and whose edge set is such that each face is connected to all the vertices along its boundary. In other words, the quadrangulation $Q$ is obtained from the map $M$ by adding a vertex
to the center of each face and then joining each such vertex to all of the vertices (counted
with multiplicity) that one encounters while tracing the boundary of that face, and then removing the edges of~$M$. Notice that
$Q$ is bipartite, with the two classes indexed by $V$ and $F$, and that all of the faces of~$Q$
are quadrilaterals, with one quadrilateral for each edge of~$M$. Let $M^*$ denote the \textit{dual map\/}
of~$M$ (the map whose vertices correspond to the faces of~$M$ --- an edge joins two vertices in
$M^*$ if an edge borders the corresponding faces in $M$).
The oriented edge $e$ of~$M$ determines an oriented edge $e_0$ of~$Q$ that has the same initial endpoint as $e$ and is the next edge clockwise
(among all edges of~$M$ and $Q$ that start at that endpoint) from $e_0$.
We refer to the endpoint of~$e_0$ in $V$ as the \textit{root\/} and the endpoint in $F$ as the \textit{dual root}.
(The dual root is the same as the outer face that we defined earlier.)

Any subset of edges $E'\subset E$ of the map $M$ defines a triangulation,
which is formed by taking the union of the quadrangulation $Q(M)$, the edges $E'$,
and the duals of the remaining edges $E\setminus E'$; each quadrilateral is subdivided into two triangles.
If each triangle is assigned unit area, then the combinatorial definition of the area of a cycle from Section~\ref{sec::intro of maps}
coincides with this geometric area.

\subsection{Peano exploration path and lattice walks}\label{sec:inside unicyle}

To further explain the hamburger-cheeseburger bijection, we briefly describe how the bijection works restricted to $\Theta^0_n$ and $\Theta^1_n$,
following \cite[Section 4.1]{sheffield2011quantum}.
We start with the bijection between $\Theta^0_n$ and $(M,e,T)$ where $T$ is a spanning tree on a rooted planar map $(M,e)\in \M_n$.
Since $T$ is a subset of the edges of~$M$ corresponding to a spanning tree, the set $T^*$ of dual edges to the edges in the complement of~$T$ is necessarily a spanning tree of~$M^*$.
The union of $T$, $T^*$, and $Q=Q(M)$ forms a triangulation, with each triangle containing two edges from $Q$ and one edge from either $T$ or $T^*$.
Let $e_0,e_1,e_2,\cdots,e_{2n}=e_0$ be the sequence of edges hit by the \textit{Peano exploration path},
which traverses each edge of~$Q$ exactly
once, with an element of~$V$ on the left and an element of~$F$ on the right, before returning to the initial edge~$e_0$.
This exploration path goes through each triangle without ever traversing an edge of~$T$ or $T^*$;
the path describes the interface between $T$ and $T^*$.

For each edge $e_i$ of the quadrangulation, let $d(e_i) = (d_1(e_i), d_2(e_i))$, where $d_1(e_i)$ is the length of the path within the tree $T$ from the $V$-endpoint of~$e_i$ to the root, and $d_2(e_i)$ is the length of the path within the dual tree $T^*$ from the
$F$-endpoint of~$e_i$ to the dual root.  Then the sequence $d(e_0),d(e_1),\cdots,d(e_{2n}) = d(e_0)$ is a
simple walk on the lattice $\Z^2_{\ge 0}$ of non-negative integer pairs, starting and ending at $(0,0)$.
To such a walk we can associate a corresponding
word in $\Theta^0_n$ by writing $\hambg$ or $\chsbg$ each time the first or second (respectively)
coordinate of~$d(e_i)$ goes up and $\hamdr$ or $\chsdr$ each time the first or second (respectively)
coordinate of~$d(e_i)$ goes down.

In the bijection, every quadrilateral of~$Q$ corresponds to a burger.
The quadrilateral is divided by an edge in $T\cup T^*$ into two triangles;
quadrilaterals of $Q$ that are divided by $T$ edges correspond to hamburgers,
while quadrilaterals divided by $T^*$ edges correspond to cheeseburgers.
The first triangle that the
path goes through corresponds to the step at which that burger was added to the stack, while
the second corresponds to the step at which the same burger was ordered.

\begin{remark}\label{rmk:edge number of tree}
Given a perfect word $W\in \Theta_n^0$ and $(M,e,T)$ which is the image of~$W$ under the bijection, from the discussion above, the number of edges in $T$ equals the number of~$\hambg$'s in $W$.  In particular, the number of vertices of $M$ is one more than the number of~$\hambg$'s in $W$.
\end{remark}

Now we describe the bijection between $\Theta^1_n$ and triples of the form $(M,e,U)$ where $U$ is a spanning unicycle on a rooted planar map $(M,e)$ of $n$ edges. Similar to the above, we define $U^*$ to be the dual edge set of $E\setminus U$ on $M^*$, which is a spanning 2-component forest of $M^*$. The interface between $U$ and $U^*$ is not a single loop, but instead consists of two loops $\L_0$ and $\L_1$, where we take $\L_0$ to be the loop that passes through the starting edge $e_0$.

Start exploring along loop $\L_0$ from the midpoint of $e_0$. $\L_0$ passes through at least one triangle that contains an edge in the loop of $U$; let $\Delta$ be the last such triangle, and let $e$ be the corresponding edge of $U$ in $\Delta$.  We let $\widetilde T=U\setminus\{e\}$ and $\widetilde T^*=U^*\cup\{e^*\}$, which ``flips the quadraleteral'', replacing the edge with the dual edge.
The effect of this replacement is to join loops $\L_0$ and $\L_1$ to form a single loop between the spanning tree $\widetilde T$ and its dual tree~$\widetilde T^*$.
To construct a word $W\in\Theta^1_n$ corresponding to $(M,e,U)$, we first apply the hamburger-cheeseburger bijection for $\Theta^0_n$ to construct a word $\widetilde W \in \Theta^0_n$ corresponding to $(M,e,\widetilde T)$ by the procedure described above.  From the construction, the quadrilateral containing $\Delta$ corresponds to a pair $\{\chsbg,\chsdr\}$. $W$ is obtained from $\widetilde W$ by replacing the $\chsdr$ in this pair by an $\frshorder$.  We call the $\chsbg$ in the pair the \textit{fresh} $\chsbg$ in $W$. By the definition of triangle~$\Delta$, the fresh $\chsbg$ will be the burger that fulfills this $\frshorder$ when reading $W$. Thus $W$ is a perfect word belonging to $\Theta^1_n$.

\enlargethispage{12pt}
 In $W$, between the fresh $\chsbg$ and the \frshorder, suppose there are $\ell$ $\hamdr$'s which are fulfilled by \hambg's before the fresh \chsbg. By the definition of perfect word and fresh order, after removing these $\hamdr$'s and those \hambg's fulfilling them, the remaining word is a perfect word of length $2n-2\ell$. Suppose there are $2m$ other letters between the fresh $\chsbg$ and the $\frshorder$ in $W$. Since the quadrilateral containing $\Delta$ is the last one to be traversed by $\L_0$ among all quadrilaterals containing edges in the loop of $U$, the length $L$ of the cycle in $U$ is $\ell+1$. Since the number of edges inside the loop is $m$, the area $A$ of the cycle in $U$ is $(1+\ell+2m)$. 
}

\section{Asymptotic enumeration of perfect words} \label{sec:edge-density}

In this section we compute the asymptotic number of perfect words in $|\Theta^k_n|$, which together with
equations \eqref{eq:E[sand]}, \eqref{eq:Var[sand]}, and \eqref{eq:Theta} gives the sandpile edge density.
In the course of characterizing~$\Theta^k_n$, we also characterize the cycles in $k$-excess graphs.
\begin{theorem}\label{thm:ell-moment}
For any fixed nonnegative integer $k$,
\begin{equation}\label{eq:ell-moment}
\lim\limits_{n\to\infty}\frac{|\Theta^k_n|}{n^k\,|\Theta^0_n|}=\frac{1}{k!\,4^k}\,.
\end{equation}
\end{theorem}

To prove Theorem~\ref{thm:ell-moment}, we study the canonical injection from $\Theta^k_n$ to $\Theta^0_n \times \binom{[2n]}{k}$:
Given a perfect word $W\in\Theta^k_n$, we can replace each of the $k$ $\frshorder$'s in $W$ with $\chsdr$'s to obtain a perfect word
in $\Theta^0_n$, and record the positions of the $\frshorder$'s as a set of $k$ distinct elements of $[2n]=\{1,\dots,2n\}$; this map is invertible so it is an injection.
A word $W$ in $\Theta^0_n$ together with a set of distinct positions $\{i_1,\dots,i_k\}$ is in the image of this injection precisely
when $W_{i_1},\dots,W_{i_k}$ are all $\chsdr$'s, and just prior to each of these orders, the top burger in the stack is a $\chsbg$.

Theorem~\ref{thm:ell-moment} can be interpreted as a statement about the probability that a random element of
$\Theta^0_n \times \binom{[2n]}{k}$ is in the image of the injection.  For a random word $W\in\Theta^0_n$ and a random position $i$,
$\P\big[W_i=\chsdr\big]=\frac14$, and at a random time, provided the stack is nonempty, the top burger on the stack is a $\chsbg$ with probability $\frac12$.  
For large $n$, it is plausible that these events at the same random time are approximately uncorrelated.
As long as $k$ is not too large ($k\ll \sqrt{n}$), we expect random
distinct positions $i_1,\dots,i_k$ to be far apart, and that
consequently these events at the times $i_1,\dots,i_k$ are nearly
independent.  Provided that this intution is correct, then
$|\Theta^k_n| \approx |\Theta^0_n| \binom{2n}{k} / 8^k$, which when
$k\ll\sqrt{n}$ would give the theorem.  In this section, we justify a more
precise version of this intuition to prove the theorem.

To make a more precise statement of this approximate independence, we consider subwords of $W\in\Theta^0_n$.
Let $W[a,b]$ denote the subword from positions $a$ through $b$ inclusive.  Let $w_j$ be the subword
\begin{equation} \label{subword}
w_j:=W[\max(i_j-s+1,1),i_j]
\end{equation}
for $j=1,\dots,k$, i.e., $w_j$ is the subword of length~$s$ which ends at position $i_j$ (unless the position is too close to the front, in which case the length will be less than $s$).
Since the perfect word $W$ corresponds to random walks in the quadrant that start and end at the origin,
we expect the subwords $w_1,\dots,w_k$ to be close in distribution to i.i.d.\ uniformly random words of length~$s$,
so long as both $k$ and $s$ are small enough for the subwords to be disjoint and not to contain enough letters
to detect that $W$ is not quite an unbiased random walk.
\begin{theorem} \label{thm:approx-independence}
  Assume $k^3 s^2 \ll n$, and consider the collection of subwords
  $(w_1,\dots,w_k)$ defined in \eqref{subword} from a random perfect
  word in $\Theta^0_n$ and independent random indices
  $i_1<\cdots<i_k$ in $\binom{[2n]}{k}$.  The
  subwords are nonoverlapping with probability $1-o(1)$ and have total
  variation distance $o(1)$ from a list of $k$ i.i.d.\ uniformly
  random words of length $s$ which are independent of the indices.
  (The $o(1)$ terms go to zero as $k^3 s^2/n\to 0$.)
\end{theorem}
\begin{proof}
  Let $i_1<\cdots<i_k$ be a random $k$-tuple uniformly drawn from
  $\binom{[2n]}{k}$, and let $\hat w_1,\ldots, \hat w_k$ be $k$ i.i.d.\
  uniformly random words in $\{\hambg,\chsbg,\hamdr,\chsdr\}^s$, which
  are also independent of the $i_1,\dots,i_k$.  Our goal is to sample
  a uniformly random perfect word in $W\in\Theta^0_n$ so that the
  subwords defined by \eqref{subword} using the indices
  $i_1,\dots,i_k$ coincide with $\hat w_1,\ldots,\hat w_k$.  Our
  strategy is to first sample an independent uniformly random perfect
  word $X\in\Theta^0_n$, and then modify it, using $i_1,\dots,i_k$ and
  $\hat w_1,\dots,\hat w_k$ and some auxillary randomness, to obtain~$W$.
  This modification proecedure defines a Markov chain, and we
  restrict ourselves to modification procedures for which the Markov
  chain satisfies detailed balance conditional on $i_1,\ldots,i_k$,
  i.e., $\P[X\to W\,|\,i_1,\ldots,i_k]=\P[W\to X\,|\,i_1,\ldots,i_k]$, so
  as to ensure that $W$ is uniformly random, and in fact uniformly random
  even conditional on $i_1,\ldots,i_k$.  After verifying detailed
  balance, we argue that with high probability $W$ has the desired subwords:
  $w_1=\hat w_1,\ldots,w_k=\hat w_k$.

  We use the following modification rule: If $i_1<s$ or
  $i_{j+1}-i_j<s$ for some $j$, then no change is made.  Otherwise,
  all the subwords have length $s$ and are disjoint.  We then start by
  overwriting the relevant positions in $X$ with the values given by
  $\hat w_1,\ldots,\hat w_k$ to obtain a new word $Y$.  The walk in $\Z^2$
  defined $Y$ is unlikely to return to its start, so some more changes
  are required to rebalance it.

  To describe this rebalancing of the walk, we
  follow Sheffield \cite{sheffield2011quantum} in using a pair of
  coordinates that are rotated $45^\circ$ from the edges of $\Z^2$:
  we let $u$ denote the number of burgers minus the number of orders,
  and $v$ denote the ``discrepancy'' between hamburgers and cheeseburgers:
  \[
  \#\{\hambg,\chsbg\}-\#\{\hamdr,\chsdr\}=u
  \quad\quad\textrm{and}\quad\quad\#\{\hambg,\chsdr\}-\#\{\hamdr,\chsbg\}=v\,.
  \]
  In these coordinates, the letters correspond to the following steps:
  \[
  \hambg=(+1,+1)\quad\quad\chsbg=(+1,-1)\quad\quad\chsdr=(-1,+1)\quad\quad\hamdr=(-1,-1)\,.
  \]
  If $u>0$, for instance, then focusing on the first coordinate, we change $u/2$ of the $+1$'s to $-1$'s,
  while ignoring the second coordinate.  The second coordinate can then be rebalanced ignoring the first
  coordinate.
  When doing this rebalancing, we only change letters
  whose position is in the range from $\frac12 n$ to $\frac32 n$, and which do not lie within the subwords.
  For reasons that will become apparent, out of all possible such ways to rebalance the walk $Y$, we pick one
  uniformly at random, and let $Z$ denote the resulting walk.  (If there are no such ways to rebalance $Y$, we
  let $W=Z=X$.)

  We will argue later that $Z$ is likely to remain within the
  quadrant, but if not, then we let $W=X$.  Otherwise, $Z$ is a
  perfect word, and we would like to take $W=Z$, but to ensure
  detailed balance conditional on $i_1,\ldots,i_k$,
  we consider $Z$ to be a proposal, and use the
  Metropolis rule to reject this proposal with some probability and
  instead take $W=X$.  The probabilities that $X$ proposes $Z$ and
  that $Z$ proposes $X$ (conditional on $i_1,\ldots,i_k$)
  are almost the same, but there is a small difference
  arising from the (likely) possibility that there are different numbers of
  ways to do the rebalancing in the two cases.  If there are $r_1$ ways to
  do the rebalancing when going from $X$ to $Z$, and $r_2$ ways when going
  from $Z$ to $X$, then the probability of accepting the proposed move is
  $\min(1,r_1/r_2)$.
  
  We have specified the process which produces the uniformly random
  perfect word $W$, and at this point we argue that with high
  probability all the steps succeed so that in the end $w_1=\hat
  w_1,\ldots,w_k=\hat w_k$.

  Since $k^2 s \ll n$, with probability $1-o(1)$ the subwords are all disjoint and have length $s$.

  After the positions are overwritten to obtain $Y$, certainly
  $|u|,|v|\leq 2 k s$.  Since $k s\ll n$, the number $m$ of letters
  with position between $\frac12 n$ and $\frac32 n$ outside the
  subwords is certainly $m=(1+o(1))n$.  For each of the two
  coordinates, standard large deviation estimates together with the
  cycle-lemma construction of random Catalan paths \cite{DM} imply that there
  are between $\frac12 m -\alpha\sqrt{m}$ and $\frac12 m
  +\alpha\sqrt{m}$ letters with $+1$ in the first coordinate (and
  similarly for the second coordinate), with probability tending to
  $1$ as $\alpha\to\infty$.  Assuming this event occurs, there is a
  way to rebalance the walk.  Provided that the rebalanced walk $Z$
  remains in the quadrant, since $|u|,|v|\leq 2 k s \ll \sqrt{m}$, the
  ratio $r_1/r_2$ tends to $1$, so the proposed move would be almost
  always accepted.

  We are left to argue that the rebalanced walk $Z$ almost always
  remains in the quadrant.  Suppose that between times $\eps n$ and
  $(2-\eps)n$ the walk $X$ remains at distance at least $h$ from the
  boundary of the quadrant.  If $\eps\ll 1/k$ and $s\leq\eps n$, then
  with probability $1-o(1)$ all of the subword regions are contained
  with the interval from $\eps n$ to $(2-\eps)n$.  There are at most
  $3 k s$ letters that get changed.  If $h\geq 6 k s$, then we would
  be guaranteed that the modified walk $Z$ would remain in the
  quadrant.

  For the initial perfect word $X$, if we eliminate the $\chsbg$ and
  $\chsdr$ letters, will be a random Catalan path of a random length
  $2\ell$ which is concentrated around $n\pm O(\sqrt{n})$.  Consider a
  uniformly random Catalan path of length $2\ell$.  As
  $\ell\to\infty$, near its endpoints it behaves as a Bessel(3)
  process.  Since a Bessel(3) process always remains positive, between
  positions $\eps\ell$ and $(2-\eps)\ell$ the Catalan path is likely
  to be at height at least $\delta\sqrt{\eps\ell}$, with probability
  tending to $1$ as $\delta\to0$.  Setting $\eps=o(1/k)$ and
  $\delta=o(1)$, we find that $X$ remains distance at least $6 k s$
  from the boundary of the quadrant with high probability provided
  $k s \ll \sqrt{n/k}$, i.e., $k^3 s^2 \ll n$.
\end{proof}

Using this approximate i.i.d.\ property of the perfect words in
$\Theta^0_n$, we can characterize $\Theta^k_n$ to prove both
Theorem~\ref{thm:ell-moment} and Theorem~\ref{thm::loop length result}.

\begin{proof}[Proof of Theorem~\ref{thm:ell-moment} and Theorem~\ref{thm::loop length result}]
So long as $k^3 s^2\ll n$, Theorem~\ref{thm:approx-independence} gives a coupling of the subwords $w_1,\dots,w_k$ of $W\in\Theta^0_n$ at random locations $i_1<\cdots<i_k\in\binom{[2n]}{k}$ with i.i.d.\ uniformly random words of length $s$, so that w.h.p.\ the subwords equal the random words.  In the i.i.d.\ words, the probability that the last letter is $\chsdr$ is $\frac14$.  As $s$ gets large, the probability that a given subword does not contain a burger that is put on the stack just prior to the $\chsdr$ is $O(1/\sqrt{s})$, and the probability that this burger is a $\chsbg$ is $\frac12$.
We can extend the i.i.d.\ words by prepending uniformly random letters, and then the probability that for each word the last letter is a $\chsdr$ that matches a $\chsbg$ is exactly $8^{-k}$, and with probability $1-O(k/\sqrt{s})$ each matching $\chsbg$ occurs within the (unextended) subword.
With probability $8^{-k}-O(k/\sqrt{s})$ the word $W$ with the inidices $(i_1,\ldots,i_k)$ is in the image of the injection,
and provided $k/\sqrt{s} \ll 8^{-k}$, this is $8^{-k}(1-o(1))$.  This proves Theorem~\ref{thm:ell-moment}.

Assuming again that $k^3 s^2\ll n$ and $k/\sqrt{s} \ll 8^{-k}$, the above coupling gives a characterization for a random perfect word $W^k$ in $\Theta^k_n$.
Suppose the fresh orders occur at positions $i_1<\cdots<i_k$, and $j_1,\ldots,j_k$ are the positions of the corresponding fresh cheeseburgers.
Then with high probability these indices alternate, $j_1<i_1<j_2<i_2<\cdots<j_k<i_k$, so that the loops given by the hamburger-cheesburger bijection are unnested, and the subwords $W^k[j_1+1,i_1-1],\ldots,W^k[j_k+1,i_k-1]$ are within $o(1)$ variation distance of $k$ i.i.d.\ simple random walks of the following type: when reading backwards, the hamburger and cheeseburger counts are always nonpositive, and the cheeseburger count ends at $0$.  These
walks are of course equivalent to walks in the quadrant, and from the discussion in Section~\ref{sec:burger bijection}, the length of the walk corresponds to the area of the cycle, and the hamburger deficit corresponds to the length of the cycle.  This proves Theorem~\ref{thm::loop length result}.
\end{proof}

\section{Sandpile density}\label{sec:small terms}

The sandpile (edge) density is a straightfoward consequence of the above lemmas:

\begin{proof}[Proof of Theorem~\ref{thm::sandpile-density} (edge density)]
  Combining equations \eqref{eq:E[sand]}, \eqref{eq:Var[sand]}, \eqref{eq:Theta}, and \eqref{eq:ell-moment}
  gives $\E^0_n[\level(\sigma(M_n))] = (1+o(1))n/4$ and $\Var^0_n[\level(\sigma(M_n))] = o(n^2)$.
  The amount of sand is $n$ more than the level, which gives the sandpile density with
  respect to the number of edges.
\end{proof}

For the random planar map with $n$ edges, the expected number of vertices is $n/2+1$.
To obtain the sandpile density with respect to the number of vertices, we need to know
that the number of vertices is sharply concentrated about its expected value.

\begin{lemma}\label{lemma::ratio of V and E}
For any fixed $k>0$, for a random planar map $M_n$ drawn from $(\M_n,\P^0_n)$,
\[\lim\limits_{n\to\infty}\E^0_n\left[\left(\frac{|E(M_n)|}{2|V(M_n)|}\right)^k\right]=1.\]
\end{lemma}
\begin{proof}
Let $J$ denote the number of $\hambg$'s contained in a uniformly random perfect word in $\Theta^0_n$.
With $\textrm{Cat}_j$ denoting the $j$th Catalan number, we have
\begin{align*}\label{eq::F1m}
|\Theta^0_n| \times \P[J=j] &= \binom{2n}{2j}\textrm{Cat}_j\textrm{Cat}_{n-j}
 = \frac{(2n)!}{j!\,(j+1)!\,(n-j)!\,(n-j+1)!} \\
 &= \binom{n+1}{j} \binom{n+1}{j+1} \frac{(2n)!}{(n+1)! (n+1)!}\,.
\end{align*}

As is well known, the binomial coefficients are sharply concentrated, with tails that are at least as small as Gaussian tails.
In particular, there are positive constants $C$ and $c$ for which
\[
\P\big[|J-\tfrac{n+1}{2}|>t\big] \leq C \, e^{-c\, t^2/n}\,.
\]

Recalling the properties of the hamburger-cheeseburger bijection
from Section~\ref{sec:burger bijection}, the number of edges in $E'$
is $|J|$, and since $E'$ forms a spanning tree, the map $M_n$ has
$|V|=|J|+1$ vertices, while of course the number of edges is $|E|=n$.
Regardless of how atypical $|J|$ may be, we have $1\leq|V|\leq n+1$, so $n/(n+1) \leq |E|/|V| \leq n$.
With, say $t=n^{2/3}$, we obtain
\begin{align*}
\E\left[(|E|/|V|)^k\right]&=\E[(|E|/|V|)^k \1_{||V|-n/2|>t}]+\E[(|E|/|V|)^k\1_{||V|-n/2|\le t}]\\
&= n^{k} e^{-\Theta(n^{1/3})}+2^k(1+O(kn^{-1/3})).\qedhere
\end{align*}
\end{proof}

\begin{proof}[Proof of Theorem~\ref{thm::sandpile-density} (vertex density)]
It suffices to prove that $\lim_{n\to\infty} \E_n^0[(\rho_v-5/2)^2] = 0$.
\begin{align*}
\rho_v-\frac{5}{2} & =2 \rho_e \left(\frac{|E|}{2|V|}-1\right)+2\left(\rho_e-\frac{5}{4}\right) \\
\left(\rho_v-\frac{5}{2}\right)^2 
 & =
 4\rho_e^2 \left(\frac{|E|}{2|V|}-1\right)^2
 + 8 \rho_e \left(\rho_e-\frac{5}{4}\right) \left(\frac{|E|}{2|V|}-1\right)
 + 4\left(\rho_e-\frac{5}{4}\right)^2
\end{align*}
Recall that $1\le \rho_e\le 2$.
By Lemma~\ref{lemma::ratio of V and E}, $\lim_{n\to \infty}\E^0_n[(|E|/|V|-2)^k]=0$ for any fixed integer $k>0$, so the first two terms above converge to $0$
in expectation.
From the edge density part of Theorem~\ref{thm::sandpile-density}, we see that the last term converges to $0$ as well.
\end{proof}

\newcommand{\arXiv}[1]{\href{http://arxiv.org/abs/#1}{arXiv:#1}}
\newcommand{\arxiv}[1]{\href{http://arxiv.org/abs/#1}{#1}}
\newcommand{\MRhref}[2]{\href{http://www.ams.org/mathscinet-getitem?mr=#1}{MR#1}}
\def\@rst #1 #2other{#1}
\newcommand\MR[1]{\relax\ifhmode\unskip\spacefactor3000 \space\fi
  \MRhref{\expandafter\@rst #1 other}{#1}}

\bibliographystyle{hmralphaabbrv}
\setlength{\itemsep}{-.3ex}


\end{document}